\documentclass[twoside]{article}%
\usepackage{amssymb}
\usepackage{amsfonts}
\usepackage{amsmath}
\usepackage{graphicx}%
\providecommand{\U}[1]{\protect\rule{.1in}{.1in}}
\newtheorem{theorem}{Theorem}

\newtheorem{remark}[theorem]{Remark}

\newenvironment{proof}[1][Proof]{\noindent\textbf{#1.} }{\ \rule{0.5em}{0.5em}}
\begin{document}

\title{\textbf{Blowup for the Euler and Euler-Poisson Equations with Repulsive
Forces}}
\author{M\textsc{anwai Yuen\thanks{E-mail address: nevetsyuen@hotmail.com }}\\\textit{Department of Applied Mathematics,}\\\textit{The Hong Kong Polytechnic University,}\\\textit{Hung Hom, Kowloon, Hong Kong}}
\date{Revised 05-Oct-2010}
\maketitle

\begin{abstract}
In this paper, we study the blowup of the $N$-dim Euler or Euler-Poisson
equations with repulsive forces, in radial symmetry. We provide a novel
integration method to show that the non-trivial classical solutions $(\rho
,V)$, with compact support in $[0,R]$, where $R>0$ is a positive constant and
in the sense which $\rho(t,r)=0$ and $V(t,r)=0$ for $r\geq R$, under the
initial condition%
\begin{equation}
H_{0}=\int_{0}^{R}rV_{0}dr>0,
\end{equation}
blow up on or before the finite time $T=R^{3}/(2H_{0})$ for pressureless
fluids or $\gamma>1.$

The main contribution of this article provides the blowup results of the Euler
$(\delta=0)$ or Euler-Poisson $(\delta=1)$ equations with repulsive forces,
and with pressure $(\gamma>1)$, as the previous blowup papers (\cite{MUK}
\cite{MP}, \cite{P} and \cite{CT}) cannot handle the systems with the pressure
term, for $C^{1}$ solutions.


\end{abstract}

\section{Introduction}

The isentropic Euler $(\delta=0)$ or Euler-Poisson $(\delta=\pm1)$ equations
can be written in the following form:
\begin{equation}
\left\{
\begin{array}
[c]{rl}%
{\normalsize \rho}_{t}{\normalsize +\nabla\cdot(\rho u)} & {\normalsize =}%
{\normalsize 0,}\\
\rho\lbrack u_{t}+(u\cdot\nabla)u]{\normalsize +\nabla}P & {\normalsize =}%
{\normalsize \rho\nabla\Phi,}\\
{\normalsize \Delta\Phi(t,x)} & {\normalsize =\delta\alpha(N)}%
{\normalsize \rho,}%
\end{array}
\right.  \label{Euler-Poisson}%
\end{equation}
where $\alpha(N)$ is a constant related to the unit ball in $R^{N}$:
$\alpha(1)=1,$ $\alpha(2)=2\pi$ and $\alpha(3)=4\pi$. And as usual, $\rho
=\rho(t,x)\geq0$ and $u=u(t,x)\in\mathbf{R}^{N}$ are the density and the
velocity respectively. $P=P(\rho)$\ is the pressure function. The $\gamma$-law
can be applied on the pressure term $P(\rho)$, i.e.%
\begin{equation}
{\normalsize P}\left(  \rho\right)  {\normalsize =K\rho}^{\gamma},
\label{gamma}%
\end{equation}
which is a common hypothesis. If the parameter\ is set as ${\normalsize K>0}$,
we call the system \textit{with pressure}; if ${\normalsize K=0}$, we call it
\textit{pressureless}. The constant $\gamma=c_{P}/c_{v}\geq1$, where $c_{P}$,
$c_{v}$\ are the specific heats per unit mass under constant pressure and
constant volume respectively, is the ratio of the specific heats, that is, the
adiabatic exponent in the equation (\ref{gamma}). In particular, the fluid is
called isothermal if $\gamma=1$. If\ ${\normalsize K>0}$, we call the system
with pressure; if ${\normalsize K=0}$, we call it pressureless.\newline In the
above systems, the self-gravitational potential field $\Phi=\Phi(t,x)$\ is
determined by the density $\rho$ itself, through the Poisson equation
(\ref{Euler-Poisson})$_{3}$.\newline When $\delta=-1$, the system can model
fluids that are self-gravitating , such as gaseous stars. In addition, the
evolution of the simple cosmology can be modelled by the dust distribution
without pressure term. This describes the stellar systems of collisionless and
gravitational $n$-body systems \cite{FT}. And the pressureless Euler-Poisson
equations can be derived from the Vlasov-Poisson-Boltzmann model with the zero
mean free path \cite{G1}. For $N=3$ and $\delta=-1$, the equations
(\ref{Euler-Poisson}) are the classical (non-relativistic) descriptions of a
galaxy in astrophysics. See \cite{BT} and \cite{C}, for details about the
systems.\newline When $\delta=1$, the system is the compressible Euler-Poisson
equations with repulsive forces. The equation (\ref{Euler-Poisson})$_{3}$ is
the Poisson equation through which the potential with repulsive forces is
determined by the density distribution of the electrons. In this case, the
system can be viewed as a semiconductor model. See \cite{Cse} and \cite{Lions}
for detailed analysis of the system.

On the other hand, the Poisson equation (\ref{Euler-Poisson})$_{3}$ can be
solved as%
\begin{equation}
{\normalsize \Phi(t,x)=\delta}\int_{R^{N}}G(x-y)\rho(t,y){\normalsize dy,}%
\end{equation}
where $G$ is Green's function for the Poisson equation in the $N$-dimensional
spaces defined by
\begin{equation}
G(x)\doteq\left\{
\begin{array}
[c]{ll}%
|x|, & N=1;\\
\log|x|, & N=2;\\
\frac{-1}{|x|^{N-2}}, & N\geq3.
\end{array}
\right.
\end{equation}

Usually, the Euler-Poisson equations can be rewritten in the scalar form:%
\begin{equation}
\left\{
\begin{array}
[c]{rl}%
\frac{\partial\rho}{\partial t}+\underset{k=1}{\overset{N}{\Sigma}}u_{k}%
\frac{\partial\rho}{\partial x_{k}}+\rho\underset{k=1}{\overset{N}{\Sigma}%
}\frac{\partial u_{k}}{\partial x_{k}} & {\normalsize =}{\normalsize 0,}\\
\rho\left(  \frac{\partial u_{i}}{\partial t}+\underset{k=1}{\overset
{N}{\Sigma}}u_{k}\frac{\partial u_{i}}{\partial x_{k}}\right)  +\frac{\partial
P}{\partial x_{i}} & {\normalsize =\rho\frac{\partial\Phi}{\partial x_{i}}%
}\text{, for }i=1,2,...N.
\end{array}
\right.  \label{gamma=1}%
\end{equation}

For the construction of the analytical solutions for the systems, interested
readers should refer to \cite{GW}, \cite{M1}, \cite{DXY}, \cite{Li} and
\cite{Y1}. The results for local existence theories can be found in \cite{M2},
\cite{B} and \cite{G}. The analysis of stabilities for the systems may be
referred to \cite{SI}, \cite{A}, \cite{E}, \cite{MUK}, \cite{MP}, \cite{P},
\cite{DLY}, \cite{DXY}, \cite{J}, \cite{Y2}, \cite{CT} and \cite{CH}.

We seek the radial symmetry solutions
\begin{equation}
\rho(t,\vec{x})=\rho(t,r)\text{ and }\vec{u}=\frac{\vec{x}}{r}V(t,r)=:\frac
{\vec{x}}{r}V\text{,}%
\end{equation}
with the radius $r=\left(  \sum_{i=1}^{N}x_{i}^{2}\right)  ^{1/2}$.

For the solutions in spherical symmetry, the Poisson equation
(\ref{Euler-Poisson})$_{3}$ is transformed to%
\begin{equation}
{\normalsize r^{N-1}\Phi}_{rr}\left(  {\normalsize t,x}\right)  +\left(
N-1\right)  r^{N-2}\Phi_{r}{\normalsize =}\alpha\left(  N\right)
\delta{\normalsize \rho r^{N-1},}%
\end{equation}%
\begin{equation}
\Phi_{r}=\frac{\alpha\left(  N\right)  \delta}{r^{N-1}}\int_{0}^{r}%
\rho(t,s)s^{N-1}ds.
\end{equation}
By standard computation, the Euler-Poisson equations in radial symmetry can be
written in the following form:%
\begin{equation}
\left\{
\begin{array}
[c]{c}%
\rho_{t}+V\rho_{r}+\rho V_{r}+\dfrac{N-1}{r}\rho V=0,\\
\rho\left(  V_{t}+VV_{r}\right)  +P_{r}(\rho)=\rho\Phi_{r}\left(  \rho\right)
.
\end{array}
\right.  \label{eq12345}%
\end{equation}

Historically, Makino, Ukai and Kawashima initially defined the tame solutions
\cite{MUK} for outside the compact of the solutions%
\begin{equation}
V_{t}+VV_{r}=0.
\end{equation}
Following this, Makino and Perthame considered the tame solutions for the
system with gravitational forces \cite{MP}. After that Perthame discovered the
blowup results for $3$-dimensional pressureless system with repulsive forces
\cite{P} $(\delta=1)$. In short, all the results above rely on the solutions
with radial symmetry:
\begin{equation}
V_{t}+VV_{r}{\normalsize =}\frac{\alpha(N)\delta}{r^{N-1}}\int_{0}^{r}%
\rho(t,s)s^{N-1}ds.
\end{equation}
And the Emden ordinary differential equations were deduced on the boundary
point of the solutions with compact support:%
\begin{equation}
\frac{D^{2}R}{Dt^{2}}=\frac{\delta M}{R^{N-1}},\text{ }R(0,R_{0})=R_{0}%
\geq0,\text{ }\dot{R}(0,R_{0})=0,
\end{equation}
where $\frac{dR}{dt}:=V$ and $M$ is the mass of the solutions, along the
characteristic curve. They showed the blowup results for the $C^{1}$ solutions
of the system (\ref{eq12345}).

Recently, Chae and Tadmor \cite{CT} showed the finite time blowup, for the
pressureless Euler-Poisson equations with attractive forces $(\delta=-1)$,
under the initial condition,%
\begin{equation}
S:=\{\left.  a\in R^{N}\right\vert \text{ }\rho_{0}(a)>0,\text{ }\Omega
_{0}(a)=0,\text{ }\nabla\cdot u(0,x(0)<0\}\neq\phi, \label{chea}%
\end{equation}
where $\Omega$ is the rescaled vorticity matrix $(\Omega_{_{0}ij})=\frac{1}%
{2}(\partial_{i}u_{0}^{j}-\partial_{j}u_{0}^{i})$ with the notation
$u=(u^{1},u^{2},....,u^{N})$ in their paper and some point $x_{0}$.

They use the analysis of spectral dynamics to show the Racatti differential
inequality,%
\begin{equation}
\frac{D\operatorname{div}u}{Dt}\leq-\frac{1}{N}(\operatorname{div}u)^{2}.
\label{ineq1}%
\end{equation}
The solution for the inequality (\ref{ineq1}) blows up on or before
$T=-N/(\nabla\cdot u(0,x_{0}(0))$.

However, their method cannot be applied to the system with repulsive forces to
obtain the similar blowup result.

On the other hand, in \cite{Y2}, we have the blowup results if the solutions
with compact support under the condition,
\begin{equation}
2\int_{\Omega(t)}(\rho\left\vert u\right\vert ^{2}+2P)dx<M^{2}-\epsilon,
\end{equation}
where $M$ is the mass of the solution.

In this article, the alternative approach is adopted to show that there is no
global existence of $C^{1}$ solutions for the system, (\ref{gamma=1})
$(\delta=0$ or $\delta=1),$ with compact support without the condition
(\ref{chea}). We notice that the conditions in our result are different from
the works of Engerlberg et. al \cite{ELT}.

\begin{theorem}
\label{thm:1 copy(1)}Consider the $N$-dimensional Euler $(\delta=0)$ or
Euler-Poisson equations with repulsive forces $(\delta=1)$
(\ref{Euler-Poisson}). The non-trivial classical solutions $\left(
\rho,V\right)  $, in radial symmetry, with compact support in $\left[
0,R\right]  $, where $R>0$ is a positive constant (which $\rho(t,r)=0$ and
$V(t,r)=0$ for $r\geq R$) and the initial velocity such that:
\begin{equation}
H_{0}=\int_{0}^{R}rV_{0}dr>0,
\end{equation}
blow up on or before the finite time $T=R^{3}/(2H_{0}),$ for pressureless
fluids $(K=0)$ or $\gamma>1$.
\end{theorem}

The solutions $(\rho,u)$ may lose their regularity, for example the velocity
function $V\in C^{0}$ only or the shock waves appear on or before the finite
time $T$.

\section{Integration Method}

In this section, we present the proof of Theorem \ref{thm:1 copy(1)}. The
technique of the proof was selected simply to deduce the partial differential
equations to the Racatti equation, to show the blowup result. However, we note
our integration method is novel to the studies of blowup for this kind of the systems.

\begin{proof}
In general, we show that the $\rho(t,x(t;x))$ preserves its positive nature as
the mass equation (\ref{gamma=1})$_{1}$ can be converted to be%
\begin{equation}
\frac{D\rho}{Dt}+\rho\nabla\cdot u=0, \label{eqq2}%
\end{equation}
with the material derivative,%
\begin{equation}
\frac{D}{Dt}=\frac{\partial}{\partial t}+\left(  u\cdot\nabla\right)  .
\label{eqq1}%
\end{equation}
We integrate the equation (\ref{eqq2})$:$%
\begin{equation}
\rho(t,x)=\rho_{0}(x_{0}(0,x_{0}))\exp\left(  -\int_{0}^{t}\nabla\cdot
u(t,x(t;0,x_{0}))dt\right)  \geq0,
\end{equation}
for $\rho_{0}(x_{0}(0,x_{0}))\geq0,$ along the characteristic curve.

We use the momentum equation (\ref{eq12345})$_{2}$ with the non-trivial
solutions in radial symmetry, $\rho_{0}\neq0$, to have:%
\begin{equation}
V_{t}+VV_{r}+K\gamma\rho^{\gamma-2}\rho_{r}=\Phi_{r},
\end{equation}%
\begin{equation}
V_{t}+\frac{\partial}{\partial r}(\frac{1}{2}V^{2})+K\gamma\rho^{\gamma-2}%
\rho_{r}=\Phi_{r},
\end{equation}%
\begin{equation}
rV_{t}+r\frac{\partial}{\partial r}(\frac{1}{2}V^{2})+K\gamma r\rho^{\gamma
-2}\rho_{r}=r\Phi_{r},
\end{equation}
with multiplying $r$ on the both sides.\newline We take integration with
respect to $r,$ to the above equation, for $\gamma>1$ or $K\geq0$:%
\begin{equation}
\int_{0}^{R}rV_{t}dr+\int_{0}^{R}r\frac{d}{dr}(\frac{1}{2}V^{2})+\int_{0}%
^{R}K\gamma r\rho^{\gamma-2}\rho_{r}dr=\int_{0}^{R}r\Phi_{r}dr,
\end{equation}%
\begin{equation}
\int_{0}^{R}rV_{t}dr+\int_{0}^{R}r\frac{d}{dr}(\frac{1}{2}V^{2})+\int_{0}%
^{R}\frac{K\gamma r}{\gamma-1}d\rho^{\gamma-1}=\int_{0}^{R}\left[
\frac{\alpha(N)\delta r}{r^{N-1}}\int_{0}^{r}\rho(t,s)s^{N-1}ds\right]  dr,
\end{equation}%
\begin{equation}
\int_{0}^{R}rV_{t}dr+\int_{0}^{R}r\frac{d}{dr}(\frac{1}{2}V^{2})+\int_{0}%
^{R}\frac{K\gamma r}{\gamma-1}d\rho^{\gamma-1}\geq0,
\end{equation}
for $\delta\geq0.$\newline It follows with integration by part:%
\begin{equation}
\int_{0}^{R}rV_{t}dr-\frac{1}{2}\int_{0}^{R}V^{2}dr+\frac{1}{2}\left[
RV(t,R)^{2}-0\cdot V(t,0)^{2}\right]  -\int_{0}^{R}\frac{K\gamma}{\gamma
-1}\rho^{\gamma-1}dr+\frac{K\gamma}{\gamma-1}\left[  R\rho^{\gamma
-1}(t,R)-0\cdot\rho^{\gamma-1}(t,0)\right]  \geq0.
\end{equation}
The above inequality with the boundary compact condition of $V(t,R)=0$ and
$\rho(t,R)=0$, becomes%
\begin{equation}
\int_{0}^{R}rV_{t}dr-\frac{1}{2}\int_{0}^{R}V^{2}dr-\int_{0}^{R}\frac{K\gamma
}{\gamma-1}\rho^{\gamma-1}dr=0.
\end{equation}
As $r$ and $t$ are independent variables and $V$ is $C^{1}$ in the domain
$[0,R]$ in the assumption of the theorem, we may change the differentiation
and the integration as the following:%
\begin{equation}
\frac{d}{dt}\int_{0}^{R}rVdr-\frac{1}{2}\int_{0}^{R}V^{2}dr-\int_{0}^{R}%
\frac{K\gamma}{\gamma-1}\rho^{\gamma-1}dr\geq0,
\end{equation}%
\begin{equation}
\frac{d}{dt}\frac{1}{2}\int_{0}^{R}Vdr^{2}-\frac{1}{2}\int_{0}^{R}\frac{1}%
{2r}V^{2}dr^{2}\geq\int_{0}^{R}\frac{K\gamma}{\gamma-1}\rho^{\gamma-1}dr\geq0,
\end{equation}
for $\gamma>1$ or $K=0.$\newline For the non-trivial initial condition
$\rho_{0}\geq0$, we have the following differential inequality:%
\begin{equation}
\frac{d}{dt}\frac{1}{2}\int_{0}^{R}Vdr^{2}-\frac{1}{2}\int_{0}^{R}\frac{1}%
{2r}V^{2}dr^{2}\geq0,
\end{equation}%
\begin{equation}
\frac{d}{dt}\int_{0}^{R}Vdr^{2}\geq\int_{0}^{R}\frac{1}{2r}V^{2}dr^{2}%
\geq\frac{1}{2R}\int_{0}^{R}V^{2}dr^{2}, \label{eq11111}%
\end{equation}%
\begin{equation}
\frac{d}{dt}\int_{0}^{R}Vdr^{2}\geq\frac{1}{2R}\int_{0}^{R}V^{2}dr^{2}.
\end{equation}
By denoting
\begin{equation}
H:=H(t)=\int_{0}^{R}rVdr=\frac{1}{2}\int_{0}^{R}Vdr^{2},
\end{equation}
and with the Cauchy-Schwarz inequality,%
\begin{equation}
\left\vert \int_{0}^{R}V\cdot1dr^{2}\right\vert \leq\left(  \int_{0}^{R}%
V^{2}dr^{2}\right)  ^{1/2}\left(  \int_{0}^{R}1dr^{2}\right)  ^{1/2},
\end{equation}%
\begin{equation}
\frac{\left\vert \int_{0}^{R}Vdr^{2}\right\vert }{R}\leq\left(  \int_{0}%
^{R}V^{2}dr^{2}\right)  ^{1/2},
\end{equation}%
\begin{equation}
\frac{4H^{2}}{R^{2}}\leq\int_{0}^{R}V^{2}dr^{2},
\end{equation}%
\begin{equation}
\frac{2H^{2}}{R^{3}}\leq\frac{1}{2R}\int_{0}^{R}V^{2}dr^{2},
\end{equation}
the inequality (\ref{eq11111}) becomes%
\begin{equation}
\frac{d}{dt}H\geq\frac{1}{2R}\int_{0}^{R}V^{2}dr^{2}\geq\frac{2H^{2}}{R^{3}},
\end{equation}%
\begin{equation}
\frac{d}{dt}H\geq\frac{2H^{2}}{R^{3}}.
\end{equation}
With the initial condition: $H_{0}=\int_{0}^{R}rV_{0}dr>0,$ we can obtain%
\begin{equation}
H\geq\frac{-R^{3}H_{0}}{2H_{0}t-R^{3}}.
\end{equation}
Therefore, the solutions blow up on or before the finite time $T=R^{3}%
/(2H_{0}).$

This completes the proof.
\end{proof}

\begin{remark}
For controlled experiments in engineering, fluids are kept in a fixed ball
solid container with a radial $R$. Therefore, it requires the compact support
condition for $t\geq0$,
\begin{equation}
\rho(t,r)=0\text{ and }V(t,r)=0\text{,}%
\end{equation}
with $r\geq R$. This corresponding condition is called no-slip condition
(solid boundary condition) \cite{Day} and \cite{CTC}.\newline On the other
hand, in computing simulations, the systems are usually coupled with the
similar boundary conditions for real applications. Therefore, the condition
for compact support (no-slip condition) is reasonable in modelings. But for
free boundary problems, fluids may not be bounded by a fixed volume for all
time. Therefore, further research is needed to study the corresponding result
in future works.
\end{remark}

\begin{remark}
It is still an open question whether or not there exists time-local $C^{1}%
$-solution with compact support for any initial condition with compact
support. On the other hand, if the global solutions with compact support whose
radii expand unboundedly as time tends to infinity, the discussion of this
paper can offer no information about this case.
\end{remark}

\begin{remark}
This article has shed new light on situations with the pressure term. In
particular, it provides the blowup results of the Euler $(\delta=0)$ or
Euler-Poisson $(\delta=1)$ equations with repulsive forces, and with pressure
$(\gamma>1)$. This is the main contribution of the article, as the previous
blowup papers (\cite{MUK} \cite{MP}, \cite{P} and \cite{CT}) cannot handle the
systems with the pressure term, for $C^{1}$ solutions. A further refinement
for the non-radial symmetry is expected in future studies.
\end{remark}

\section{Acknowledgement}

The author would like to thank the comments of reviewers to improve the article.

\end{document}